\documentclass[twocolumn]{autart}    

\makeatletter
\def\@biblabel#1{}
\makeatother

\usepackage[fleqn]{amsmath}
\usepackage{graphicx,mathrsfs,amssymb,cite,enumerate,mathrsfs,balance,float}
\usepackage{color}
\usepackage[round]{natbib}

\usepackage{subfigure}

\begin{document}

\begin{frontmatter}
\runtitle{Novel Multi Agent Models for Chemical Self-assembly}  

\title{Novel Multi-Agent Models for Chemical Self-assembly \thanksref{footnoteinfo}} 

\thanks[footnoteinfo]{the National Natural Science Foundation
of China under grants 91427304, 61673373 and 11688101, the National
Key Basic Research Program of China (973 program) under grant
2014CB845301/2/3, and the Leading research projects of Chinese Academy
of Sciences under grant QYZDJ-SSW-JSC003.}

\author[Chen,Ning]{Zheng Ning}\ead{ningzheng15@mails.ucas.ac.cn},    
\author[Chen]{Ge Chen}\ead{chenge@amss.ac.cn},
\author[Bullo]{Francesco Bullo}\ead{bullo@engineering.ucsb.edu},

\address[Chen]{National Center for Mathematics and Interdisciplinary Sciences \& Key Laboratory of Systems and
Control, Academy of Mathematics and Systems Science, Chinese Academy of Sciences, Beijing 100190,
China}
\address[Ning]{School of Mathematical Sciences, University of
Chinese Academy of Sciences, Beijing 100049, China}  
\address[Bullo]{Department of Mechanical Engineering and the Center of Control, Dynamical-Systems and Computation, University of California at Santa Barbara,
CA 93106-5070, USA}


\begin{abstract}                          
The chemical self-assembly has been considered as one of most important scientific problems in the 21th Century; however, since the process of self-assembly is very complex,
there is few mathematic theory for it currently. This paper provides a novel multi-agent model for chemical self-assembly, where the interaction between agents adopts the classic Lennard-Jones potential.
Under this model, we propose an optimal problem by taking the temperature as the control input, and choosing the internal energy as the optimal object. A numerical solution
for our optimal problem is also developed. Simulations show that our control scheme can improve the product of self-assembly. Further more, we give a strict analysis for the
self-assembly model without noise, which corresponds to an attraction-repulsion multi-agent system, and prove it converges to a stable configuration eventually.
\end{abstract}

\end{frontmatter}

\section{Introduction}

Self-assembly is the process in which disordered components form an organized structure with local interactions among the components without external forces \citep{Whitesides2002Self}. It reveals how disordered components form the ordered structure in nature and the understanding of self-assembly could help us create nano-structured materials and build new nanostructures.
This topic has raised a lot of interest in physics, chemistry and biology in recent several decades. In the \emph{Science} 125th anniversary, the magazine raised 125 important scientific problems in 21st century \citep{Service2005How}. Among these problems, they picked 25 most important ones, one of which is ``How Far Can We Push Chemical Self-Assembly?"

Recent developments of self-assembly in chemistry have been made in both experiments  \citep{Nykypanchuk2008DNA} and computational simulations \citep{Wilber2009Self}. Experimental scientists have
carried on a large number of self-assembly experiments,  including molecular, nanoparticles and protein molecular.
A lot of assembly products have been discovered. Theoretical scientists use computers to simulate the self-assembly process\citep{klotsa2013controlling}. However, since
there are too many kinds of assembly products, and the assembly process is very complex, the hidden laws of self-assembly are hard to find through experiments and simulations.
Currently, there is few mathematic theory for the self-assembly. This paper tries to model the self-assembly by multi-agent systems and build an optimal control on them.

Multi-agent systems composed by multiple interacting agents have drawn considerable attention from various fields in the past two decades. In physics, the synchronization phenomena of coupled oscillators, flashing fireflies, and chirping crickets is investigated \citep{Kuramoto1975International,acebron2005kuramoto}; in biology, scientists model animal flocking behavior \citep{buhl2006disorder,vicsek1995novel}; in sociology, the emergence and spread of public opinions can also be investigated as multi-agent system\citep{deffuant2000mixing,hegselmann2002opinion}. A central issue of multi-agent system study is to understand how local interactions among the elements lead to collective behavior of the whole group. Because of the importance, effort has been devoted to the mathematical analysis of collective behavior of multi-agent systems \citep{chen2017critical}.
 A new method called as ``soft control" has been proposed which keeps the local rule of the existing agents in the system and controls the collective behavior indirectly by changing exterior environment \citep{han2006soft}.

In this paper, we model the dynamics of self-assembly as a multi-agent system, in which each agent denotes a component, and the interaction between agents denotes the force between components. Because the real self-assembly is very complex, we need make some simplification. In our model, we assume the agents are homogeneous and isotropic, and the interaction between agents adopts the classical Lennard-Jones potential. Such a system corresponds to some practical systems such as the self-assembly of gold nanoparticles\citep{gittins2002dense}. Because the temperature is crucial in self-assembly, we treat it as a control input, and try to find its optimal control to the assembly product we want. This control method can also be considered as a soft control, however different from the method of adding some special agents in previous work \citep{han2006soft}.

The contribution of this paper can be formulated as follows:
First, in this paper we propose a new mathematical model for chemical self-assembly.
As mentioned above, the mathematical models and analysis are very important to understand the behind law of chemical self-assembly,
however they are very few currently.
Our model provides an enter point for the analysis of self-assembly. Although our model is ideal, it still keeps the essential feature of chemical self-assembly.
Also, it is possible that we could extend the model to some real assembly experiments like the self-assembly of gold nanoparticles \citep{kaplan2006intermolecular}.

Secondly, we put forward a mathematical method to find the optimal temperature control for our model. It is well-known that the temperature is a critical value in chemical systems.
How to find the optimal temperature control is an important issue in the research of chemical self-assembly. Currently chemists mainly adopt empirical methods, and still lack the guidance of mathematical theory. Our method is a way to get the optimal control for our model, and hope it can be extended to some real chemical systems.

Finally, we give a strict analysis for the multi-agent self-assembly model without noise, and prove it will converge to a stable configuration eventually.
The noise-free self-assembly model can be treated as a multi-agent attractive-repulsive model, which has attracted a lot of interest in the study of flocking algorithms.
Our system and results can provide some new idea for the modeling and analysis to the flocking research.

The rest of this article is organized as follows: In Section \ref{sec_model} we introduce a multi-agent model for chemical self-assembly and give a result for the noise-free case. In Section \ref{sec ctrl} we propose an optimal control problem to our model and explore a numerical solution.  Section \ref{sec nf} provide some simulations using our control laws, while Section \ref{sec cln} concludes this paper.

\section{A multi-agent model for self-assembly}\label{sec_model}

This section proposes a novel multi-agent model for the chemical self-assembly. The model considers $N$ homogeneous and isotropic particles move in a fluid,
and each particle $i$ contains a position variable $X_i=X_i(t)\in \mathbb{R}^{3}$ and a velocity variable
$V_{i}=V_i(t)\in \mathbb{R}^{3}$.  We use the classic Langevin equation to formulate the dynamics of the particles,
 which is,
for $t\geq 0$ and $1\leq i\leq N$,  the position and velocity of particle $i$ is driven by
\begin{eqnarray}\label{noise model}
\left\{%
\begin{array}{l}
\dot{X}_{i}=V_{i}\\
\dot{V}_{i}=-BV_{i}+f_{i}+\xi_{i}
\end{array}%
\right.,
\end{eqnarray}
where $B$ is a constant denoting the damping coefficient of each particle in the fluid, $f_{i}=f_i(t)$ is the force of particle $i$ affected by other particles, and $\xi_{i}=\xi_i(t)$ denotes the Brownian force produced by the thermal noise.

In practical chemical self-assembly, the interaction between particles is very complex (e.g. Van der Waals, capillary, $\pi-\pi$ hydrogen bonds).
 To be simplified, this paper assumes the interaction between particles is additive, and can be described by the classic  Lennard-Jones (L-J) potential.
In mathematics, the L-J potential between two particles $i$ and $j$ can be formulated as
\begin{eqnarray}\label{LJpotential}
\Phi_{ij}=\Phi_{ij}(t)=\varepsilon\Big(\frac{r_{m}}{r_{ij}}\Big)^{12}-2\varepsilon\Big(\frac{r_{m}}{r_{ij}}\Big)^6,
\end{eqnarray}
where $\varepsilon$ and $r_m$ are  constants denoting the depth of the potential well and the distance at which the potential reaches its minimum respectively,
and
$$r_{ij}=r_{ij}(t)=\|X_i-X_j\|_2$$ denotes the distance between the particle $i$ and $j$ at time $t$. Here $\|\cdot\|_2$ denotes the Euclidean norm.  Thus, the force of  particle $i$ affected by particle $j$
is
\begin{eqnarray*}
\begin{aligned}
f_{ij}=f_{ij}(t)=-\nabla\Phi_{ij}&=12 \varepsilon \left( \frac{r_{m}^{6}}{r_{ij}^{7}} -   \frac{r_{m}^{12}}{r_{ij}^{13}} \right) \left(-\nabla r_{ij} \right)\\
&=12\varepsilon \left( \frac{r_{m}^{6}}{r_{ij}^{8}} -   \frac{r_{m}^{12}}{r_{ij}^{14}} \right)  \left( X_{j}-X_{i} \right),
\end{aligned}
\end{eqnarray*}
where $\nabla g$ denotes  the gradient of  function $g$ in $\mathbb{R}^{3}$ (i.e., $\nabla g=(\frac{\partial g}{\partial x},\frac{\partial g}{\partial y}, \frac{\partial g}{\partial z})$).
Consider the forces between particles are additive, so the total force of particle $i$ affected by other particles is
\begin{eqnarray}\label{force}
\begin{aligned}
f_{i}&=f_i(t)=\sum_{j \neq i}f_{ij}\\
&=\sum_{j \neq i} 12 \varepsilon \left( \frac{r_{m}^{6}}{r_{ij}^{8}} -   \frac{r_{m}^{12}}{r_{ij}^{14}} \right)  \left( X_{j}-X_{i} \right).
\end{aligned}
\end{eqnarray}

According to the Langevin equation theory \citep{Coffey2004The}, the Brownian force $\xi_{i}$ has zero mean and its
 covariance is
 \begin{eqnarray*}
\mbox{Cov}(\xi_i(t),\xi_j(t'))
=\left\{%
\begin{array}{ll}
0,  & \mbox{if~} i\neq j\\  %
2Bk_{b}U\delta(t-t'), & \mbox{otherwise}
\end{array}%
\right.,
\end{eqnarray*}
where $k_{b}$ denotes the Boltzmann constant, $U$ denotes the absolute temperature, and
$\delta(t)$ denotes the Dirac delta function.

The system (\ref{noise model}) seems ideal, however it keeps the essential feature of chemical self-assembly. Also, it is possible to extend this model to some real assembly systems like the self-assembly of gold nanoparticles.

For the convenience of mathematical analysis, the system (\ref{noise model}) can be transformed to a group of stochastic differential equations.
By the Langevin equation theory \citep{Coffey2004The}, the integration of the force $\xi_i$, $\int_{0}^{t}\xi_i(\tau)\mathrm{d}\tau$, has the same property as
$\sqrt{2Bk_{b}U} W_{i}(t)$, where $W_i(t)$ is a standard Wiener process independent with $\{W_j(t)\}_{j\neq i}$.
Therefore, $\xi_i(t) \mathrm{d}t= \mathrm{d}[\int_{0}^{t}\xi_i(\tau)\mathrm{d}\tau]$ can be written as $\sqrt{2Bk_{b}U}\mathrm{d} W_{i}(t)$.
and the system (\ref{noise model}) can be written as the following stochastic differential equations:
\begin{eqnarray}\label{nm_sde}
\left\{%
\begin{array}{l}
\mathrm{d}{X}_{i}={V}_{i}\mathrm{d}t\\
\mathrm{d}{V}_{i}=(-BV_{i}+f_{i})\mathrm{d}t+\sqrt{2Bk_{b}U}\mathrm{d} W_{i}(t)
\end{array}%
\right..
\end{eqnarray}

\section{Optimal control for self-assembly model}\label{sec ctrl}

The Hamiltonian plays a key role in a physical system.
In this paper our prime goal is to minimize the Hamiltonian, which indicates the assembly product reaches a most stable state.
By (\ref{LJpotential}), the Hamiltonian of our system (\ref{noise model}) is
\begin{eqnarray}\label{H}
H=H(t)=\frac{1}{2}\sum_{i=1}^{N} V_{i}^{T}V_{i}+\sum_{i<j} \varepsilon \left( \frac{r_{m}^{12}}{r_{ij}^{12}}-2\frac{r_{m}^{6}}{r_{ij}^{6}} \right).
\end{eqnarray}

We also need to choose a suitable control input. In real chemical self-assembly, the particles are very small, and are hard to be controlled directed; however we can control the
external environment to intervene the assembly product. Temperature control is very important in chemical self-assembly. At present, it is mainly based on empirical method and lacks mathematical theory guidance. In this paper, we try to build a method on how to control the temperature to minimize the Hamiltonian.
This method can  be also treated as a soft control which initially proposed by \citep{han2006soft}.

To be simplified we set $u(t):=k_{b}U(t)$ as the control input, and rewrite the system (\ref{nm_sde}) as the following dynamics:
\begin{eqnarray}\label{nmu sde}
\left\{%
\begin{array}{l}
\mathrm{d}{X}_{i}={V}_{i}\mathrm{dt}\\  %
\mathrm{d}{V}_{i}=(-BV_{i}+f_{i})\mathrm{d}t+\sqrt{2Bu}\mathrm{d}W_i(t)
\end{array}%
\right..
\end{eqnarray}
We aim to minimize the  Hamiltonian $H(T)$ with $T$ being a fixed time. According to (\ref{H}) and (\ref{nm_sde}), $H(t)$ is a stochastic process and we calculate the differential of $H(t)$.
\begin{eqnarray}\label{dh00}
\begin{aligned}
\mathrm{d} H=&\frac{1}{2}\sum_{i=1}^{N} \mathrm{d} ( V_{i}^{T}V_{i}) + \frac{1}{2}\varepsilon\sum_{i \neq j}  \mathrm{d}  \left( \frac{r_{m}^{12}}{r_{ij}^{12}}-2\frac{r_{m}^{6}}{r_{ij}^{6}} \right).
\end{aligned}
\end{eqnarray}
The first part $\frac{1}{2}\sum_{i=1}^{N} \mathrm{d} ( V_{i}^{T}V_{i})$ involves $V_i(t)$, we must use It\^{o}'s formula to calculate its differential. According to the general It\^{o}'s formula
(Theorem 4.2.1 in \citep{Oksendal1985Stochastic}),
\begin{eqnarray}\label{dh11}
\begin{aligned}
&\frac{1}{2}\sum_{i=1}^{N} \mathrm{d} ( V_{i}^{T}V_{i})\\
&=\frac{1}{2}\sum_{i=1}^{N} ( 2 V_{i}^{T}\mathrm{d}V_{i}+\mathrm{d}V_{i}^{T}\mathrm{d}V_{i} )   \\
&= \frac{1}{2} \sum_{i=1}^{N}\left[ 2\left( -BV_{i}+f_{i} \right)^{T}V_{i} \mathrm{d}t+2\sqrt{2Bu}V_i \mathrm{d}W_i(t)  \right] \\
&~~ +\frac{1}{2}\sum_{i=1}^{N}\mathrm{d}V_{i} ^{T}\mathrm{d}V_{i}  \\
&=\sum_{i=1}^{N}(-BV_{i}^{T}V_{i})\mathrm{d}t+f_{i}^{T}V_i \mathrm{d}t +\sqrt{2Bu}V_i \mathrm{d}W_i(t)\\
&~~ +\frac{1}{2}\mathrm{d}V_{i} ^{T}\mathrm{d}V_{i} .
\end{aligned}
\end{eqnarray}
According to the stochastic differential equation theory, we have $\mathrm{d}t \cdot \mathrm {d}t=\mathrm{d}t \cdot \mathrm{d}W(t)=\mathrm{d}W(t) \cdot \mathrm {d}t =0$, and $\mathrm{d}W(t) \cdot \mathrm{d}W(t) =\mathrm{d}t$. Then,
$\mathrm{d}V_{i} ^{T}\mathrm{d}V_{i}$ can be expanded into follows:
\begin{eqnarray}\label{dh12}
\begin{aligned}
&\mathrm{d}V_{i} ^{T}\mathrm{d}V_{i} \\
&= \left[\left( -BV_{i}^{T}+ f_{i}^{T}\right) \mathrm{d}t+\sqrt{2Bu} \mathrm{d}W_i(t) \right] ^{T}\\
   & ~~~~\times \left[\left( -BV_{i}^{T}+ f_{i}^{T}\right) \mathrm{d}t+\sqrt{2Bu} \mathrm{d}W_i(t) \right] \\
&=2Bu\mathrm{d}W_i(t)\mathrm{d}W_i(t) \\
&= 2Bu\mathrm{d}t.
\end{aligned}
\end{eqnarray}
The second part in (\ref{dh00}) only involves $X_i(t)$, It\^{o}'s formula is not needed.
\begin{eqnarray}\label{dh21}
\begin{aligned}
&\frac{1}{2}\varepsilon\sum_{i \neq j}  \mathrm{d}  \left( \frac{r_{m}^{12}}{r_{ij}^{12}}-2\frac{r_{m}^{6}}{r_{ij}^{6}} \right)=6\sum_{i \neq j} \varepsilon \left(  \frac{r_{m}^{6}}{r_{ij}^{7}}  -  \frac{r_{m}^{12}}{r_{ij}^{13}} \right)  \mathrm{d}r_{ij}  \\
&=6\sum_{i \neq j} \varepsilon  \left(  \frac{r_{m}^{6}}{r_{ij}^{7}} -  \frac{r_{m}^{12}}{r_{ij}^{13}} \right)\frac{1}{r_{ij}}\left( X_{i}-X{j} \right)^{T}\left (V_{i}-V_{j} \right) \mathrm{d}t \\
&=6\sum_{i \neq j} \varepsilon  \left(  \frac{r_{m}^{6}}{r_{ij}^{8}} -  \frac{r_{m}^{12}}{r_{ij}^{14}} \right)  \\
&~~~~\times \left[  \left(X_{i}-X{j} \right)^{T} V_{i} +\left(X_{j}-X{i} \right)^{T}V_{j} \right] \mathrm{d}t \\
&= \left[ \frac{1}{2}\sum_{i=1}^{N}\left(-f_{i}^{T}V_{i}\right)+\frac{1}{2}\sum_{j=1}^{N}\left(-f_{j}V_{j} \right) \right] \mathrm{d}t\\
&=\sum_{i=1}^{N}-f_{i}^{T}V_{i}\mathrm{d}t.
\end{aligned}
\end{eqnarray}
By (\ref{dh11}), (\ref{dh12}) and (\ref{dh21}), the differential of $H(t)$ is
\begin{eqnarray}\label{dh3}
\begin{aligned}
\mathrm{d} H=&\sum_{i=1}^{N}(-BV_{i}^{T}V_{i} +Bu )\mathrm{d}t+\sqrt{2Bu}V_i\mathrm{d}W_{i}(t).
\end{aligned}
\end{eqnarray}
The Hamiltonian  at time $T$ could be represented as follows:
\begin{eqnarray}\label{H_T}
\begin{aligned}
H(T)=& H(0)+\int_{0}^{T} \mathrm{d}H \\
=&H(0)+\int_{0}^{T} \sum_{i=1}^{N}(-BV_{i}^{T}V_{i} +Bu )\mathrm{d}t \\
&+\int_{0}^{T}\sum_{i=1}^{N}\sqrt{2Bu}V_i\mathrm{d}W_{i}(t) .
\end{aligned}
\end{eqnarray}
 Because the Hamiltonian  $H$ is stochastic, we use $E[H(T)]$ as the optimization object. Because the expectation of It\^{o}'s integral is zero, from (\ref{H_T}) we have
\begin{eqnarray} \label{expectation}
\mathrm{E} H(T)=\mathrm{E}H(0)+\int_{0}^{T}B\sum_{i=1}^{N}(-E V_{i}^{T}V_{i} +u )\mathrm{d}t.
\end{eqnarray}
Our aim is to find the optimal control $u(t)$ to minimize the value of $\mathrm{E}H(T)$.

Because the temperature is limited by an allowable region in a chemical experiment, we assume the lower and upper bounds of $u$ are $u_{\min}$ and  $u_{\max}$ respectively. Then, we consider the following optimization problem:
\begin{eqnarray}\label{opt 1}
\begin{aligned}
 \min ~& \mathrm{E}H(0)+\int_{0}^{T} B\sum_{i=1}^{N}(-E V_{i}^{T}V_{i} +u )\mathrm{d}t \\
 s.t. &\  u_{\min} \leq u \leq u_{\max}, \\
& \mathrm{d}{V}_{i}=(-BV_{i}+f_{i})\mathrm{d}t+\sqrt{2Bu}\mathrm{d}W_i(t),  i=1,...,N.
\end{aligned}
\end{eqnarray}
By (\ref{expectation}) the optimization objective $\mathrm{E}H(T)$ is a nonlinear function,
and by (\ref{nmu sde}) the velocity $V_i(t)$ is a very complex stochastic process which depends on not only  the control input $u(t)$, but also the states of other particles. It is hard to get the analytic optimal solution. As an alternative, we will develop a numerical method to optimize $\mathrm{E}H(T)$.

\subsection{Numerical method for the optimal control problem (\ref{opt 1})}

Firstly we use Monte Carlo method to transform  the stochastic constraints in (\ref{opt 1}) into deterministic  constraints.
Set $X(t):=(X_1(t),\ldots,X_N(t))$, $V(t):=(V_1(t),\ldots,V_N(t))$, and $W(t):=(W_1(t),\ldots,W_N(t))$. Let $W(0:T)$ denote the trajectory $W(t),0\leq t\leq T$.
We randomly select $M$ sample trajectories
$\{X^k(0),V^k(0), W^k(0:T)\}_{1\leq k\leq M}$
 from the sample space of $\{X(0), V(0), W(0:T)\}$,
 and  change the stochastic constraints in (\ref{opt 1}) into the following constraints:
 \begin{eqnarray}\label{opt 2}
\begin{aligned}
\mathrm{d}{V}_{i}^k&=(-B V_{i}^k+f_{i}) \mathrm{d}t+\sqrt{2Bu}\mathrm{d}W_i^k, \\
&~\forall 1\leq i\leq N, 1\leq k\leq M,
\end{aligned}
\end{eqnarray}
 where $V_i^k$ and $W_i^k$ denote the $i$-th element of $V^k$ and $W^k$ corresponding to  particle $i$.

If $M$ is large enough, the objective function in (\ref{opt 1}) can be approximated by
\begin{eqnarray} \label{expectation2}
\begin{aligned}
&\mathrm{E}H(0)+\int_{0}^{T} B\sum_{i=1}^{N}(-E V_{i}^{T}V_{i} +u )\mathrm{d}t  \\
&\approx \frac{1}{M} \sum_{k=1}^{M} \Big[H(0)+ B\int_{0}^{T}\sum_{i=1}^{N}(- (V_{i}^k)^{T}V_{i}^k +u )\mathrm{d}t\Big].
\end{aligned}
\end{eqnarray}

Secondly, we  discretize the time interval $[0,T]$ into $N_T+1$ points $0,\Delta t,\ldots N_T \Delta t$ with
$\Delta t:=\frac{T}{N_T}$. Let $t_n:=n \Delta t$.
Then, we use ($V_{i}^{k}(0),V_{i}^{k}(\Delta t),\ldots,V_{i}^{k}(N_{T} \Delta t)$) to approximate the trajectory of $V_i^k(0:T)$,
and $(u(0),u(\Delta t),...,u(N_T \Delta t))$ to approximate the continuous-time control $u(0:T)$.  Notice that
\begin{eqnarray*}
X_i(t)=\int_{0}^{t}V_i(s)\mathrm{d}s+X_i(0),
\end{eqnarray*}
so at time $t_n$, the position of each particle $i$ can be approximated by
\begin{eqnarray}\label{pos_app}
X_i(t_n)\approx \Delta t \sum_{q=1}^{n}V_i(t_q)+X_i(0).
\end{eqnarray}
By (\ref{force}) and (\ref{pos_app}), the force $f_i(t_n)$ can also be approximated by
\begin{eqnarray}\label{for_app}
f_i(t_n) \approx f_{i}\Big(\Delta t \sum_{q=1}^{n}V_{1}(t_q),\ldots,\Delta t\sum_{q=1}^{n}V_{N}(t_q)\Big).
\end{eqnarray}
Corresponding, the differential $d W_i^k(t_n)$ is approximated by the difference
$W^{k}_i(t_{n})-W^{k}_i(t_{n-1})$, where the difference are independent random variables which have normal distribution and its variance is $\Delta t$ \citep{Glasserman2004Monte}. From this and (\ref{for_app}), the differential equation
 (\ref{opt 2}) can be approximated by the following difference equation:
\begin{eqnarray}\label{opt 3}
\begin{aligned}
&V_{i}^{k}(t_{n+1})-V_{i}^{k}(t_n)\approx -BV_{i}^{k}(t_{n+1}) \Delta t \\
&+f_{i}(\Delta t\sum_{q=1}^{n}V_{1}^k(t_q),\ldots, \Delta t\sum_{q=1}^{n}V^k_{N}(t_q))\Delta t \\
&+\sqrt{2Bu(t_{n+}1)}(W^{k}_i(t_{n+1})-W^{k}_i(t_n)),\\
& \    \ i=1,\ldots,N, k=1,\ldots,M,n=0,\ldots,N_{T}-1.
\end{aligned}
\end{eqnarray}
Similarly, we discretize the right side of (\ref{expectation2}) and get
\begin{eqnarray}\label{target}
\begin{aligned}
&\mathrm{E}H(0)+\int_{0}^{T} B\sum_{i=1}^{N}(-E V_{i}^{T}V_{i} +u )\mathrm{d}t  \\
&\approx \frac{1}{M}\sum_{k=1}^{M}\big[ \mathrm{E}H(0)  -B\Delta t \sum_{n=1}^{N_{T}} \sum_{i=1}^{N}V_{i}^{k}(t_n)^{T}V_{i}^{k}(t_n)\\
&~~~~+NB\Delta t  \sum_{n=1}^{N_{T}}u(t_n)\big].
\end{aligned}
\end{eqnarray}
By (\ref{target}) and (\ref{opt 3}), the optimal control problem  (\ref{opt 1}) can be approximated by the follows:
  \begin{eqnarray}\label{opt reform}
\begin{aligned}
 \min &~\frac{1}{M}\sum_{k=1}^{M}\big[ \mathrm{E}H(0)-B\Delta t\sum_{n=1}^{N_{T}} \sum_{i=1}^{N}V_{i}^{k}(t_n)^{T}V_{i}^{k}(t_n) \\
 &+NB\Delta t  \sum_{n=1}^{N_{T}}u(t_n)\big]. \\
 s.t. &\ u_{\min} \leq u(t_n) \leq u_{\max}, \\
&V_{i}^{k}(t_{n+1})-V_{i}^{k}(t_n)=-BV_{i}^{k}(t_{n+1}) \Delta t \\
&+f_{i}(\Delta t\sum_{q=1}^{n}V_{1}^k(t_q),\ldots, \Delta t\sum_{q=1}^{n}V^k_{N}(t_q))\Delta t \\
&+\sqrt{2Bu(t_{n+1})}(W^{k}_i(t_{n+1})-W^{k}_i(t_n)),\\
& \    \ i=1,\ldots,N, k=1,\ldots,M,n=0,\ldots,N_{T}-1.
\end{aligned}
\end{eqnarray}
We can solve (\ref{opt reform}) by sequential quadratic programming (SQP) method which is suitable for non-linear optimization with constraints \citep{Spellucci1998An}.

\subsection{A numerical solution and comparison with natural cooling}\label{Simulations}

\begin{figure}[ht]
  \centering
  \includegraphics[width=2.5in]{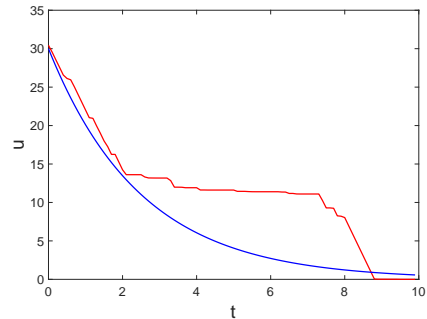}\\
  \caption{The optimal control for N=30, M=100 and t $\in [0,10]$.  }\label{noise1}
\end{figure}

In this section, we provide a numerical example to solve the problem (\ref{opt reform}).
Choose the particle number $N=30$, the stop time $T=10$, and the time step $\Delta t=0.1$.
Set $\varepsilon$ and $r_m$ in (\ref{force}) to be 3 and 2 respectively. Let the damping coefficient  $B=2$.
The initial position $X_i(0)(1\leq i\leq 30)$ is assumed to be uniformly and independently distributed in $[0,10]^3$.
For the initial velocity
$V_i(0) (1\leq i\leq 30)$, we assume $V_{ij}(0) (1\leq j\leq 3)$ has independent normal distribution whose expectation is zero and variance is $4$. Choose the lower bound $u_{\min}$ and the upper bound $u_{\max}$ of the control $u$ to be $0$ and $50$ respectively.
To effectively solve the problem (\ref{opt reform}), we adopt the annealing temperature control which assumes the temperature is non-increasing.  A solution to the problem (\ref{opt reform}) is shown by the blue curve in Figure \ref{noise1}.

\begin{figure}[ht]
  \centering
  \includegraphics[width=2.5in]{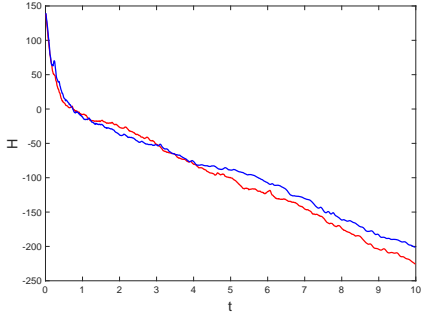}\\
  \caption{The comparison between our solution and natural cooling. The red line denotes the Hamiltonian curve using the control of our solution, while the blue one is the Hamiltonian curve
 using the natural cooling control.}\label{noise2}
\end{figure}

\begin{figure}[htbp]
\begin{center}
\subfigure[The system state at $t=10$ under the control of natural cooling.]{
\includegraphics[scale=0.28]{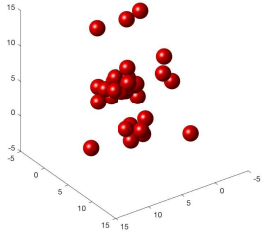}}
\subfigure[The system state at $t=10$ under the control of our numerical solution.]{
\includegraphics[scale=0.28]{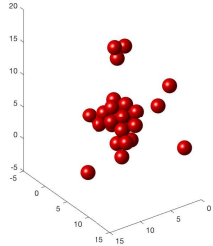}}
 \caption{Simulations of the system (\ref{noise model}) with $N=30$.}{\label{Figcom}}
\end{center}
\end{figure}

We also compare our solution with the natural cooling which is the traditional annealing control in real chemical experiments. The temperature curve of natural cooling can be approximated as the well-known Newton's law of cooling. The curve according to Newton's law of cooling is shown by the red curve in Fig. \ref{noise1}. With the same initial configuration, we run the system (\ref{noise model}) under the control of both our numerical solution and the natural cooling, where the Hamiltonian curves and system states are shown in
Figs. \ref{noise2} and \ref{Figcom}. These simulations show that our numerical solution has better performance than the natural cooling.

\section{Convergence of noise-free self-assembly model}\label{sec nf}
\renewcommand{\thesection}{\arabic{section}}
In this section we consider the noise-free case of our system (\ref{noise model}), which can be treated as a multi-agent attractive-repulsive model interested by the study of flocking algorithms \citep{Reynolds1987Flocks,Olfati2006Flocking,Cucker2011A}. Remark that the interaction between particles in our system is different from the previous works.
From (\ref{noise model}), the noise-free self-assembly model can be formulated as follows:
\begin{eqnarray}\label{model}
\left\{%
\begin{array}{l}
\dot{X}_{i}=V_{i}\\
\dot{V}_{i}=-BV_{i}+f_{i}
\end{array}%
\right., ~~1\leq i\leq N.
\end{eqnarray}

We give a convergence result for the system (\ref{model}):
\begin{thm}[Convergence of noise-free model]\label{thm}
Consider the noise-free self-assembly model  (\ref{model}).
For any initial state $(X(0),V(0))$, $(X(t),V(t))$ converges to an equilibrium point $(X(\infty),V(\infty))$ which satisfies $V_i(\infty)=\mathbf{0}$ and $f_i(\infty)=\mathbf{0}$ for all $1\leq i\leq n$.
\end{thm}

\noindent\textbf{Proof.}
The derivative of the Hamiltonian
\begin{eqnarray}\label{thm_1}
\begin{aligned}
\dot{H}&=\frac{1}{2} \sum_{i=1}^{N} \left( \dot{V_{i}^{T}}V_{i} + V_{i}^{T}\dot{V_{i}} \right) +6\sum_{i \neq j} \varepsilon \left(  \frac{r_{m}^{6}}{r_{ij}^{7}}  -  \frac{r_{m}^{12}}{r_{ij}^{13}} \right)  \dot{r_{ij}}\\
&=\sum_{i =1}^{N} -BV_{i}^{T}V_{i}+\frac{1}{2}\sum_{i=1}^{N}\left( f_{i}^{T}V_{i}+V_{i}^{T}f_{i} \right) \\
& \  \ +6\sum_{i \neq j} \varepsilon \left(  \frac{r_{m}^{6}}{r_{ij}^{8}}  -  \frac{r_{m}^{12}}{r_{ij}^{14}} \right) \left( X_{i}-X{j} \right)^{T}\left (V_{i}-V_{j} \right)\\
&=\sum_{i =1}^{N} -BV_{i}^{T}V_{i}+\frac{1}{2}\sum_{i=1}^{N}\left( f_{i}^{T}V_{i}+V_{i}^{T}f_{i} \right) \\
& \ \ +\frac{1}{2}\sum_{i=1}^{N}\left(-f_{i}^{T}V_{i}\right)+\frac{1}{2}\sum_{j=1}^{N}\left(-f_{j}V_{j} \right)\\
&=\sum_{i =1}^{N} -BV_{i}^{T}V_{i} \leq 0.
\end{aligned}
\end{eqnarray}

So the Hamiltonian $H$ will not decrease for any initial state.
According to the LaSalle invariance principle, the system will reach the state satisfying $\dot{H}=0$, which is the same as $V_{i}=\mathbf{0}$, so for any initial state we have
\begin{eqnarray}\label{thm_2}
\lim_{t \to \infty}V_{i}(t)=\mathbf{0}.
\end{eqnarray}

It remains to prove $\lim_{t \to \infty}f_{i}(t)=\mathbf{0}$ for any $1\leq i\leq N$. We prove this result by contradiction.
If there exists $i\in\{1,\ldots,N\}$ such that $f_{i}(t)$ does not converge to zero, then there exists a constant $\varepsilon>0$, an integer $k\in\{1,2,3\}$, and an infinite
sequence $t_1<t_2<\cdots $ satisfying $t_{j+1}-t_j\geq 1$ and
\begin{eqnarray}\label{thm_3}
| f_{ik} (t_j) | \ge 2\varepsilon.
\end{eqnarray}
By (\ref{thm_2}), there exists an integer $J>0$ such that
\begin{eqnarray}\label{thm_4}
| V_{ik} (t_j)| \leq \frac {1}{B} \varepsilon,~~~~\forall j\geq J.
\end{eqnarray}
Substituting this into (\ref{model}) we can get
\begin{eqnarray}\label{thm_5}
| \dot V_{ik} (t_j)|=| -B V_{ik}(t_j)+f_{ik}(t_j)| \geq \varepsilon,~~~~\forall j\geq J.
\end{eqnarray}
For convenience we omit $t$ in the next part. On the other hand,
\begin{eqnarray}\label{thm_8}
\begin{aligned}
\left\| \ddot V_{i} \right\| &=\left\|  -B \dot V_{i}+\dot f_{i} \right\| \\
&=\left\|  -B(-BV_{i}+f_{i})+\dot f_{i} \right\| \\
&=\left\|  B^{2}V_{i}-Bf_{i}+\dot f_{i} \right\| \\
&\leq \left\| B^{2}V_{i} \right\|+\left\| B f_{i} \right\|+\left\| \dot f_{i} \right\|
\end{aligned}
\end{eqnarray}
Because $\lim_{t \to \infty}V_{i}(t)=\mathbf{0}$ , the first part $\left\|B^{2}V_{i}\right\|$ is uniformly bounded. For the second part,
\begin{eqnarray}\label{thm_9}
\begin{aligned}
\left\| B f_{i} \right\|&=B \left\| \sum_{j \neq i} 12 \varepsilon  \left(  \frac{r_{m}^{6}}{r_{ij}^{8}} -  \frac{r_{m}^{12}}{r_{ij}^{14}} \right)  \left( X_{j}-X_{i} \right) \right\|  \\
&=12 B \varepsilon \sum_{j \neq i}  \left |    \frac{r_{m}^{6}}{r_{ij}^{8}} - \frac{r_{m}^{12}}{r_{ij}^{14}} \right | \left\| \left( X_{j}-X_{i} \right) \right\| \\
&=12 B \varepsilon \sum_{j \neq i}  \left|   \frac{r_{m}^{6}}{r_{ij}^{8}} -  \frac{r_{m}^{12}}{r_{ij}^{14}} \right | r_{ij}  \\
& \leq 12B \varepsilon N \max_{i\neq j} \left |   \frac{r_{m}^{6}}{r_{ij}^{7}} -  \frac{r_{m}^{12}}{r_{ij}^{13}}  \right |  \\
\end{aligned}
\end{eqnarray}
The limitation of this function of $r_{ij}$ :$\left | \frac{r_{m}^{6}}{r_{ij}^{7}}-  \frac{r_{m}^{12}}{r_{ij}^{13}} \right | $ is zero when $r_{ij} \to \infty$. So it is bounded when $r_{ij} $ is large.For another, if $r_{ij}$ is small, we can deduce it has a lower bound using the Hamiltonian $H$.
Because $H$ is decreasing. so it will always be smaller than the initial Hamiltonian $H(0)$. And also the two-particle potential has the minimum $- \varepsilon$.$r_{ij}$ will always be subject to
\begin{eqnarray}\label{thm_10}
\varepsilon \left( \frac{r_{m}^{12}}{r_{ij}^{12}}-2\frac{r_{m}^{6}}{r_{ij}^{6}} \right) +\frac{N(N-1)}{2}(-\varepsilon) \leq H(0)
\end{eqnarray}
So the $r_{ij}(t)$ has a positive low bounded, which indicates that the second part $|B f_{i}(t)|$ of the last line of (\ref{thm_8}) has a uniform upper bound.

For the third part $\dot f_{i}(t)$ of the last line of (\ref{thm_8}),according to  (\ref{force}), is
\begin{eqnarray}\label{thm_11}
\begin{aligned}
\left\| \dot f_{i}(t)\right\|&=12 \varepsilon \sum_{j \neq i}  \bigg\|  \left(  \frac{r_{m}^{6}}{r_{ij}^{8}} - \frac{r_{m}^{12}}{r_{ij}^{14}} \right)  \left( V_{j}-V_{i} \right)  \\
& +\left( - 8 \frac{r_{m}^{6}}{r_{ij}^{9}} +14  \frac{r_{m}^{12}}{r_{ij}^{15}} \right ) \frac{1}{r_{ij}} \left( X_{j}-X_{i} \right)^{T} \\
& ~~ \left( V_{j}-V_{i} \right) \left( X_{j}-X_{i} \right) \bigg\| \\
& \leq 12  \varepsilon  \sum_{j \neq i} \Bigg| \left(  \frac{r_{m}^{6}}{r_{ij}^{8}} - \frac{r_{m}^{12}}{r_{ij}^{14}} \right)  2 \sup_{i,t} ||V_{i}|| \Bigg| \\
& +\bigg| -8 \frac{r_{m}^{6}}{r_{ij}^{9}} +14  \frac{r_{m}^{12}}{r_{ij}^{15}} \bigg | \cdot 2 |r_{ij}| \sup_{i,t} ||V_{i}||.
\end{aligned}
\end{eqnarray}
Because from (\ref{thm_2}) we get $V_i$ has a uniform upper bound, and from (\ref{thm_10}) we get $r_{ij}$ has a positive low bounded,
so by (\ref{thm_11}) we obtain  $|\dot f_{i}|$ has a uniform upper bound. Given the discussion above we get $|\ddot V_{i} (t)|$ is uniformly bounded.

Since $\dot V_{ik}(t)$ is a derivable function, and $\ddot V_{ik}(t)$ is uniformly bounded,
we can find a constant $\delta>0$ such that
\begin{eqnarray}\label{thm_6}
| \dot V_{ik} (t)-\dot V_{ik} (t_j)|\leq \frac{\varepsilon}{2},~~~\forall j\geq J, t\in [t_j-\delta,t_j+\delta].
\end{eqnarray}
By (\ref{thm_6}) and (\ref{thm_5}) we get
\begin{eqnarray}\label{thm_7}
\begin{aligned}
&|V_{ik}(t_j+\delta)-V_{ik}(t_j-\delta)|\\
&=\Big|\int_{t_j-\delta}^{t_j+\delta}  \dot V_{ik} (t) dt\Big|\geq \delta  \varepsilon, ~~\forall j\geq J,
\end{aligned}
\end{eqnarray}
which is contradictory with (\ref{thm_2}). \hfill $\Box$

\begin{figure}[htbp]
\begin{center}
\subfigure[The initial state]{\label{Init}
\includegraphics[scale=0.28]{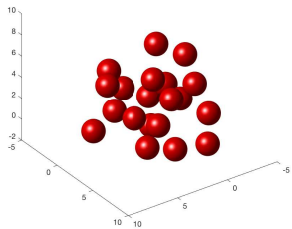}}
\subfigure[The final state]{\label{Fina}
\includegraphics[scale=0.28]{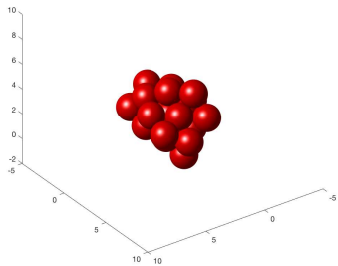}}
 \caption{A simulation for the system (\ref{model}) with $N=20$.}{\label{Fig1}}
\end{center}
\end{figure}

We use $20$ particles to simulate the system (\ref{model}) as follows: Assume the initial positions of all particles are uniformly and independently distributed in $[0,10]^{3}$,  and the initial velocities  are uniformly and independently distributed in $[0,1]^{3}$. Let $B=1$, $\varepsilon=1$, and $r_{m}=2$.  The initial state and final state of the system (\ref{model}) are shown in
Figure \ref{Fig1}.

\begin{figure}[htbp]
\begin{center}
\includegraphics[scale=0.28]{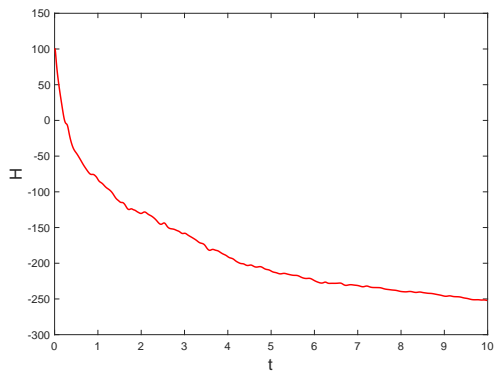}
 \caption{The Hamiltonian $H$ of the system (\ref{model}) with the same configuration as Fig \ref{Fig1}.}{\label{Fig3}}
\end{center}
\end{figure}

%
%
\section{Conclusion and future works}\label{sec cln}

This paper provides a novel multi-agent model for chemical self-assembly. Because the particles in our model cannot be controlled directly, we propose a optimal  temperature  control problem
which can be treated as a kind of soft control. A numerical method to our optimization problem is explored. Also, we consider the noise-free case of our system and prove its convergence.

The temperature control plays a key role in chemical self-assembly however few mathematical theory exists. In the future we could adjust our model and method according to the real chemical experiment. For example, the force  $f_i$ could be remodeled according to the dynamics of the assembly particles, and the damping  coefficient $B$ could be determined by the kind of the fluid.
Another future work could change the objective function in our model according to the  requirement of the real chemical experiment. For example, if we want the components assemble a specific structure, we need to choose  a suitable objective function to fit the target structure.


\bibliographystyle{plainnat}

\bibliography{bibfile}

\end{document}